\documentclass{amsart}

\usepackage{amsmath}
\usepackage{amssymb}
\usepackage{eucal}
\usepackage{amsthm}
\usepackage{verbatim}
\usepackage[hyphens]{url}
\usepackage{natbib}

\usepackage{graphicx}
\usepackage{subfigure}

\newtheorem{theorem}{Theorem}[section]

\newtheorem{lemma}[theorem]{Lemma}
\newtheorem{corollary}[theorem]{Corollary}

{\theoremstyle{definition}

}

{\theoremstyle{remark}
\newtheorem{remark}[theorem]{Remark}
}

\DeclareMathOperator{\Prob}  {\mathsf{P}}
\DeclareMathOperator{\Expec} {\mathsf{E}}

\newcommand{\mfM}{\mathfrak{M}}

\newcommand{\mcC}{\mathcal{C}}

\newcommand{\mbR}{\mathbb{R}}
\newcommand{\mbS}{\mathbb{S}}

\newcommand{\stfrac}[2]{ \leavevmode\ensuremath{
    \kern.1em\raise.5ex\hbox{\footnotesize #1}
    \kern-.1em/
    \kern-.15em\lower.25ex\hbox{\footnotesize #2}} }
\newcommand{\smfrac}[2]{ \leavevmode\ensuremath{
    \kern.1em\raise.5ex\hbox{\footnotesize $#1$}
    \kern-.1em/
    \kern-.15em\lower.25ex\hbox{\footnotesize $#2$}} }
\newcommand{\Stfrac}[2]{ \leavevmode\ensuremath{
    \kern.1em\raise.5ex\hbox{#1}
    \kern-.1em/
    \kern-.15em\lower.25ex\hbox{#2}} }
\newcommand{\Smfrac}[2]{ \leavevmode\ensuremath{
    \kern.1em\raise.5ex\hbox{$#1$}
    \kern-.25em/
    \kern-.15em\lower.5ex\hbox{$#2$}} }

\begin{document} 


\title[Linear Latent Structure Analysis]
      {Linear Latent Structure Analysis:
       Mixture Distribution Models with Linear Constraints}

\author[M.Kovtun]{Mikhail~Kovtun}
\address{Duke University, CDS\\
         2117 Campus Dr., Durham, NC, 27708}
\email{mkovtun@cds.duke.edu}

\author[I.Akushevich]{Igor~Akushevich}
\address{Duke University, CDS\\
         2117 Campus Dr., Durham, NC, 27708}
\email{aku@cds.duke.edu}

\author[K.G.Manton]{Kenneth~G.~Manton}
\address{Duke University, CDS\\
         2117 Campus Dr., Durham, NC, 27708}
\email{kgm@cds.duke.edu}

\author[H.D.Tolley]{H.~Dennis~Tolley}
\address{Brigham Young University, Department of Statistics\\
          206 TMCB, Provo, UT, 84602}
\email{tolley@byu.edu}

\thanks{This work was supported in part by NIH grants
        P01-AG17937-05 and R01-AG01159-28.}

\thanks{This article was submitted for publication
        in {\em Statistical Methodology}.}


\begin{abstract} 
A new method for analyzing high-dimensional categorical data,
Linear Latent Structure (LLS) analysis, is presented.
LLS models belong to the family of latent structure models,
which are mixture distribution models constrained to satisfy
the local independence assumption.
LLS analysis explicitly considers a family of mixed distributions as 
a linear space and LLS models are obtained by imposing linear
constraints on the mixing distribution.

LLS models are identifiable under modest conditions and are
consistently estimable.
A remarkable feature of LLS analysis is the existence of a
high-performance numerical algorithm, which reduces 
parameter estimation to a sequence of linear algebra problems.
Preliminary simulation experiments with a prototype of the algorithm
demonstrated a good quality of restoration of model parameters.
\end{abstract} 

\maketitle 

\section{Introduction} 
\label{sec:Introduction}

We present a new statistical method,
Linear Latent Structure (LLS) analysis,
which belongs to a domain of latent structure analysis.
In presentation of the method and in investigation of its
properties we follow a new understanding of latent structure
models as mixture distribution models
\citep{Bartholomew:2002}.

The latent structure analysis considers a number of categorical
variables measured on each individual in a sample,
and it is aimed to discover properties of a population
as well as properties of individuals composing a population.
The main assumption of the latent structure analysis is the
{\em local independence} assumption. 
Being formulated in a more contemporary way, this means that 
the observed joint distribution of categorical random variables 
is a mixture of independent distributions 
(see section \ref{sec:LLSModels} for more detail).
The mixing distribution is considered as a distribution
of {\em latent variable(s)}, which is thought of as containing hidden
information regarding the phenomenon under consideration.
The goal of the latent structure analysis is to discover
properties of latent variables; different approaches to this
problem are described in
\cite{Lazarsfeld:1950a,Lazarsfeld:1950b,Lazarsfeld:1968,
Goodman:1978,Langeheine:1988,Clogg:1995,Heinen:1996,Bartholomew:1999,
Marcoulides:2002}.

The various branches of latent structure analysis differ
in additional assumptions regarding latent variables---or,
equivalently, regarding mixing distribution.
Latent class analysis assumes that the mixing distribution is
concentrated in a finite number of points (called ``latent classes'').
Latent trait analysis (LTA) tries to represent mixed distributions
as a function of a latent trait, which in most cases is assumed
to be one-dimensional parameter.
In the 1990s, multidimensional latent traits were
investigated \citep{Hoijtink:1997, Reckase:1997}.
However, the application of multidimensional LTA is not as
broad as one-dimensional LTA, since estimating parameters
in multidimensional case requires additional assumptions.

The novelty of our approach lies in consideration of the space
of mixed distributions as a linear space. This allows us to
employ geometric intuition and clearly formulate the main
additional assumption of LLS that the mixed distribution is supported
by a linear subspace of the space of independent distributions.

Models arising from this approach
are identifiable under modest conditions and are
consistently estimable.
Further, there exists a high-performance
numerical algorithm, which reduces estimation of model parameters
to a sequence of linear algebra problems.
Preliminary simulation experiments with a prototype of the algorithm
(presented in section \ref{sec:SimulationStudies})
demonstrated a good quality of restoration of model parameters.

The word ``linear'' in the name of the method reflects at least
three aspects of our method:
first, the model is obtained by imposing linear constraints
on the mixing distribution;
second, the algorithm for model construction uses methods of
linear algebra;
and third, the most interesting ideas of our approach arise
from consideration of a space of distributions as a linear space.

Historically, the predecessor of LLS analysis was
Grade of Membership (GoM) analysis, which was introduced in 
\cite{Woodbury:1974}; see also \cite{Manton:1994} for detailed
exposition and additional references.
Our work on LLS analysis originated from attempts to find
conditions for consistency of GoM estimators.
The development eventually lead to a new class of models,
which differ from GoM models in a way how the model is
formulated, methods of model estimation, meaning of estimators
and their interpretation.
We present here this new class of models under the name 
``linear latent structure analysis.''


The present article concentrates on the exposition of main ideas 
of the LLS analysis and investigation of its statistical properties.
Section \ref{sec:LLSModels} describes the basics of LLS analysis.
A new procedure for parameter estimation
is explained in Section \ref{sec:EstimationOfLLSModel}.
Sections \ref{sec:IdentifiabilityOfLLS} and \ref{sec:ConsistencyOfLLS}
answer the questions concerning identifiability of LLS models
and consistency of the estimators.
Section \ref{sec:SimulationStudies} gives results of preliminary
experiments with a prototype of the algorithm implemented
by the authors.
The article is concluded by section \ref{sec:Discussion},
where we discuss some interesting properties of LLS analysis
and compare it with other kinds of latent structure analysis.


\section{LLS models} 
\label{sec:LLSModels}

The input to LLS analysis is outcome of $J$ categorical
measurements, each made on $N$ individuals.

Mathematically, we consider $J$ categorical random variables
$X_1, \dots, X_J$. The set of possible values of random
variable $X_j$ is $\{ 1, \dots L_j \}$.
This structure is described by an integer vector
$L = (L_1, \dots, L_J)$.
Two numbers, which are used frequently in the rest of
the article, are associated with that structure:
$|L| = L_1 + \dotsc + L_J$ and $|L^*| = L_1 \cdot \dotsc \cdot L_J$.

To denote response patterns, we use integer vectors
$\ell = ( \ell_1, \dots, \ell_J )$. The $j^\text{th}$ component,
$\ell_j$, represents the value of random variable $X_j$
(thus, $\ell_j \in [1..L_j]$). Note that there are $|L^*|$
different vectors $\ell$.

The joint distribution of random variables $X_1, \dots, X_J$
is given by $|L^*|$ elementary probabilities

\begin{equation}   
\label{eq:ElemProbs}
p_\ell = \Prob ( X_1 = \ell_1 \text{ and } \dotsc 
                              \text{ and } X_J = \ell_J )
\end{equation}   

In general, no restrictions other than 
$p_\ell \ge 0$ for all $\ell$ and $\sum_\ell p_\ell = 1$
are imposed on the family of elementary probabilities;
thus, one needs $|L^*|-1$ parameters to describe a joint
distribution of $X_1, \dots, X_J$.

Note that probabilities $p_\ell$ are directly estimable
from the observed data: frequencies $f_\ell = \frac{N_\ell}{N}$
(where $N$ is the total number of individuals in the sample
and $N_\ell$ is the number of individuals who responded with
response pattern $\ell$) are consistent and efficient estimators
for $p_\ell$.

Among all joint distributions one can distinguish 
{\em independent} distributions, i.e. distributions,
in which random variables $X_1, \dots, X_J$ are mutually
independent. This means that for every set of indices
$j_1, \dots, j_p$ and for every response pattern $\ell$ the relation

\begin{equation}   
\label{eq:Independence1}
\Prob ( X_{j_1} = \ell_{j_1} \text{ and } \dotsc 
                             \text{ and } X_{j_p} = \ell_{j_p} )
= \Prob( X_{j_1} = \ell_{j_1} ) \cdot \ldots \cdot
    \Prob( X_{j_p} = \ell_{j_p} )
\end{equation}   

\noindent
holds. Equation (\ref{eq:Independence1}) allows us to describe
an independent distribution using fewer parameters. 
Namely, let $\beta_{jl} = \Prob(X_j=l)$.
Then for every response pattern $\ell$,

\begin{equation}   
\label{eq:IndParam1}
p_\ell = \prod_{j=1}^J \beta_{j \ell_j}
\end{equation}   

Thus, every independent distribution can be identified with a
point $\beta = (\beta_{jl})_{jl} \in \mbR^{|L|}$. Not every
point $\beta \in \mbR^{|L|}$ corresponds to a probability distribution;
to describe a distribution, $\beta$ must satisfy the conditions:

\begin{equation}   
\label{eq:BetaCond}
\begin{cases}
    \sum_{l=1}^{L_j} \beta_{jl} = 1     & \text{for every $j$} \\
    \beta _{jl} \ge 0                   & \text{for every $j$ and $l$}
\end{cases}
\end{equation}   

Conditions (\ref{eq:BetaCond}) define a convex
$(|L|-J)$-dimensional polyhedron in $\mbR^{|L|}$, which we
denote $\mbS^L$.

Now we are ready to formulate the first assumption of the LLS
analysis:

\begin{enumerate}
\item[(G1)] \emph{
    The observed distribution is a mixture of independent distributions,
    i.e. there exist a probabilistic measure $\mu_\beta$, supported
    by $\mbS^L$, such that for every response pattern $\ell$
    \begin{equation}   
    \label{eq:LocalInd1}
        p_\ell = \int 
            \left( {\textstyle \prod_{j=1}^J \beta_{j \ell_j} } \right) 
                                                \, \mu_\beta(d\beta)
    \end{equation}   
}
\end{enumerate}

\begin{remark}
Here and later, we use measures instead of probability density
functions or cumulative distribution functions, 
as this allows us to avoid discussion of possible singularities
and specifics of the space on which probability distribution is
defined. 
This is just a matter of convenience; the above integral
can be written as
$\int \left( \prod_{j=1}^J \beta_{j \ell_j} \right) \, dF(\beta)$,
where $F(\beta)$ is the cumulative distribution function of
the mixing distribution.
\end{remark}

Assumption (G1) is a cornerstone of latent structure analysis
(often, it is called the {\em local independence} assumption).
There are a lot of excellent books and articles devoted to 
latent structure analysis; we refer to 
\cite{Lazarsfeld:1950a,Lazarsfeld:1950b,Lazarsfeld:1968,
Goodman:1978,Langeheine:1988,Clogg:1995,Heinen:1996,Bartholomew:1999,
Marcoulides:2002}
for discussion of the meaning and applicability of this assumption.

Mathematically, the assumption (G1) alone does not imply much.
For almost each distribution $(p_\ell)_\ell$ there exist
infinitely many mixing distributions $\mu_\beta$ that produce
the same observed distribution.
Thus, one needs more assumptions to make the model identifiable.
In latent structure analysis, such assumptions are usually
formulated in the form of restrictions on the support of the mixing
distribution $\mu_\beta$.

The specific assumption of LLS analysis is:

\begin{enumerate}
\item[(G2)] \emph{
    The mixing distribution $\mu_\beta$ is supported by
    a linear subspace $Q$ of $\mbR^{|L|}$.
}
\end{enumerate}

For comparison, the corresponding assumption of latent class analysis
is that $\mu_\beta$ is supported by a finite number of points
(latent classes).

When dimensionality of $Q$ is sufficiently smaller
than $|L|$, LLS model is almost surely identifiable and consistently
estimable from data. It is discussed in subsequent sections.

Informally, the existence of low-dimensional support of measure 
$\mu_\beta$ means that all measurements reflect the same underlying
hidden entity. In \cite{Kovtun:2005} we have shown that
the existence of low-dimensional support is equivalent to the existence
of a $K$-dimensional random vector $G$ such that regressions
of all indicator random vectors $Y_j$ on $G$ are linear.
($Y_j = (Y_{j1},\dots,Y_{jL_j})$; $Y_{jl}=1$, if $X_j=l$; otherwise,
$Y_{jl}=0$.)

Distributions satisfying the condition (G2) may be expected when random
variables $X_1,\dots,X_J$ represent responses to survey or exam
questions. Here, questions are intentionally chosen to discover
a single (potentially multidimensional) quantity---like 
``quality of life'' or ``mathematical knowledge''.
This is the natural domain of applications of latent structure
analysis in general, and LLS analysis in particular.

We say that a distribution is generated by a $K$-dimensional
LLS model, if it can be represented as a mixed distribution
satisfying (G2) with $\dim(Q)=K$.


\section{Estimation of LLS model} 
\label{sec:EstimationOfLLSModel}

``To define a LLS model'' means to define mixing distribution 
$\mu_\beta$, which, in turn, means specifying the supporting
subspace $Q$ and the distribution over it.

The supporting subspace may be consistently estimated from 
the observed data, i.e. the estimated subspace converges to the true one
when sample size tends to infinity. The identifiability conditions are
rather straightforward: if the dimensionality of supporting subspace is
of order of $\left(\frac{|L|-J}{2}-\max L_j \right)$ or smaller, the
supporting subspace is almost surely identifiable (theorem
\ref{th:Ident}).

LLS analysis uses nonparametric approach to description of
the mixing distribution. Thus, the knowledge about the mixing 
distribution is expressed in the form of a family of conditional 
moments of order up to $J$. 
Using these moments, the mixing distribution may be
approximated as an empirical distribution. The examples given
in section \ref{sec:SimulationStudies} demonstrate the goodness
of such approximation; see also section \ref{sec:Discussion}
for discussion of properties of this approximation.

The technical details of what follows are given in \cite{Kovtun:2005}.
Here we formulate the most important facts and pay more attention
to the most significant statistical properties such as
identifiability of the model and consistency of the estimates.

Let $K$ be the dimensionality of $Q$ and let
$\lambda^1 = (\lambda^1_{jl})_{jl}, \dots, 
\lambda^K = (\lambda^K_{jl})_{jl}$
be a basis of $Q$.
Let $g = (g_1,\dots,g_K)$ be coordinates of points of $Q$
written in the basis $\lambda^1,\dots,\lambda^K$.
This means that for points contained in $Q$
coordinates $\beta$ and $g$ connected as:

\begin{equation}   
\label{eq:BetaG}
\beta_{jl} = \textstyle{ \sum_{k=1}^K \lambda^k_{jl} \cdot g_k }
\end{equation}   

\noindent
or, in matrix form, $\beta = \Lambda g$, where $\Lambda$
is $|L| \times K$ matrix, $\Lambda = (\lambda^k_{jl})^k_{jl}$.

Recall that the support of mixing measure $\mu_\beta$ is also
restricted to $\mbS^L$ (a polyhedron defined by conditions
(\ref{eq:BetaCond})), i.e. $\mu_\beta$ is supported by intersection
of $Q$ and $\mbS^L$.
We consider only bases in which all $\lambda^k$ belong to $\mbS^L$.
In this case, coordinates $g$ of points belonging to the support
of $\mu_\beta$ satisfy $g_1+\dots+g_K=1$;
thus, $g$ are homogeneous coordinates of points from $Q \cap \mbS^L$.
It is possible to exclude any coordinate $g_1,\dots,g_K$ and
use the remaining $K-1$ coordinates to denote points of
$Q \cap \mbS^L$; however, we prefer to use the redundant set of
coordinates to preserve symmetry of equations.

Let $\mu_g$ be the measure $\mu_\beta$ written in coordinates $g$.
This means that for every function $\phi$ defined on $Q$ one has
$\int \phi(\beta) \,\mu_\beta(d\beta)
= \int \phi(\Lambda g) \,\mu_g(dg)$.
In particular,

\begin{equation}   
\label{eq:IndParam2}
p_\ell = \int 
    \left( {\textstyle \prod_{j=1}^J \beta_{j \ell_j} } \right) 
                                                \, \mu_\beta(d\beta)
= \int
    \left( {\textstyle \prod_{j=1}^J
                \sum_{k=1}^K \lambda^k_{jl} \cdot g_k } \right)
                                                \, \mu_g(dg)
\end{equation}   

Every probabilistic measure on $n$-dimensional euclidean space 
may be considered as a distribution law of an $n$-dimensional random
vector. Let $B=(B_{jl})_{jl}$ be a random vector corresponding 
to measure $\mu_\beta$ and let $G=(G_k)_k$ be a random vector 
corresponding to measure $\mu_g$.
In fact, $B$ and $G$ are the same random vector, but written
in different coordinates.

It might be shown that $X_1,\dots,X_J$ and $G$ (or $B$) have
a joint distribution; thus, one can speak about conditional 
probabilities and conditional expectations.

Some moments of order $J$ of $B$ coincide with elementary
probabilities (\ref{eq:ElemProbs}). 
Name\-ly, due to (\ref{eq:LocalInd1}),

\begin{equation}   
\label{eq:BMoments1}
M_\ell(B) = \int 
       \left( {\textstyle \prod_{j=1}^J \beta_{j \ell_j} } \right) 
                                                \, \mu_\beta(d\beta)
= p_\ell
\end{equation}   

The above equation may be extended to moments of order lower
than $J$. To proceed, we need to extend $\ell$-notation.
From now, we allow 0's in some positions of vector $\ell$. Such
0's mean that we ``do not care'' about values of corresponding
random variables. A vector $\ell$ with some components equal to 0
may be also thought of as a set of all response patterns,
which have arbitrary values on ``do not care'' places
and coincide with $\ell$ on all other places.
Then $p_\ell$ will be a marginal probability,
and the corresponding moments of $B$ will be:

\begin{equation}   
\label{eq:BMoments2}
M_\ell(B) = \int 
       \bigg( \prod_{j \,:\, \ell_j \neq 0} \beta_{j \ell_j}  \bigg) 
                                                \, \mu_\beta(d\beta)
= p_\ell
\end{equation}   

The set of moments $M_\ell(B)$ for all $\ell$ (including
$\ell$ with zeros) is all what is directly estimable from
the observation.

Another set of values of interest is the set of conditional
moments of order $v=(v_1,\dots,v_K)$ of random vector $G$:

\begin{equation}   
\label{eq:GikDef}
g^v_\ell =
    \Expec ( G_1^{v_1} \cdot \ldots \cdot G_K^{v_K} \mid X = \ell )
\end{equation}   

\noindent
Here $\Expec$ denotes the expectation, 
and $X=\ell$ is an abbreviation for conjunction of conditions
$X_j = \ell_j$ for all $j$ such that $\ell_j \neq 0$.
Note that the values $g^v_\ell$ depend on the choice 
of the basis $\lambda^1,\dots,\lambda^K$.

These conditional moments express the knowledge regarding individuals
that can be obtained from the measurements.
In particular, conditional expectations (equation (\ref{eq:GikDef1})
below) may be considered as estimators of individual coordinates
in state space (see also section \ref{sec:Discussion}).

Among all conditional moments, moments of order 1, or conditional
expectations, have special importance, and we use special notation
for them:

\begin{align}   
\label{eq:GikDef1}
&g_{\ell 1} \overset{\text{\scriptsize{def}}}{=}
    g^{(1,0,\dots,0)}_\ell = \Expec(G_1 \mid X=\ell) \notag \\
&\dots \\
&g_{\ell K} \overset{\text{\scriptsize{def}}}{=}
    g^{(0,\dots,0,1)}_\ell = \Expec(G_K \mid X=\ell) \notag
\end{align}   

The above values satisfy the following equation 
\citep[][section 6.2]{Kovtun:2005}:

\begin{equation}   
\label{eq:MainEq}
M_\ell(B) \cdot \left(
    \lambda^1_{jl} \cdot g^{v^1}_\ell + \dots
    + \lambda^K_{jl} \cdot g^{v^K}_\ell
\right)
=
M_{\ell'}(B) \cdot g^v_{\ell'}
\end{equation}   

\noindent
Here: (a) $v^k$ denotes a vector $v$ with $k^\text{th}$ component 
increased by 1 (for example, if $v=(1,3,2,1)$, then $v^3=(1,3,3,1)$),
(b) response pattern $\ell$ must have 0 at $j^\text{th}$ position,
(c) $\ell'$ denotes the response pattern obtained from $\ell$ by
replacing 0 at $j^\text{th}$ position by $l$
(for example, if $\ell = (1,0,0,2,1)$ and $j=3$, then 
$\ell' = (1,0,l,2,1)$).
Equation (\ref{eq:MainEq}) holds for every $j$, $l$, $v$,
and every $\ell$ containing 0 at $j^\text{th}$ place.

A special case of equation (\ref{eq:MainEq}) when $v=(0,\dots,0)$
is:

\begin{equation}   
\label{eq:MainEqExp}
M_\ell(B) \cdot \left(
    \lambda^1_{jl} \cdot g_{\ell 1} + \dots
    + \lambda^K_{jl} \cdot g_{\ell K}
\right)
=
M_{\ell'}(B)
\end{equation}   

\noindent
The right-hand side of this equation does not involve $g^v_\ell$
because $g^{(0,\dots,0)}_\ell = 1$.

By combining all equations (\ref{eq:MainEq}) for all possible
$v$ and $\ell$ with normalization equations
(like $\sum_k g_{\ell k} = 1$) one obtains
the main system of equations \citep[see][section 7]{Kovtun:2005}.

The important property of the main system of equations
is given by the following

\begin{theorem}   
\label{th:MainSys1}
Let $M_\ell(B)$ be moments of a distribution generated by
$K$-dimensio\-nal LLS model.
Let $\lambda^1,\dots,\lambda^K$ be any basis of the supporting
subspace $Q$ and let $g^v_\ell$ be conditional moments
calculated with respect to this basis.

Then $\lambda^k$ and $g^v_\ell$ give a solution of the main
system of equation with coefficients $M_\ell(B)$.
\end{theorem}   

Moreover, for almost all (in the strict mathematical sense) 
distributions {\em every} solution of the main system of equations 
is a basis of the supporting subspace and conditional moments,
calculated with respect to this basis. This implies that LLS
model is almost surely identifiable; we discuss this fact in more
detail in the next section.

Note that equations (\ref{eq:MainEq}) are linear with respect to
variables $g^v_\ell$. Thus, if one knows a basis of $Q$,
it is sufficient to solve a linear system of equations to find
the conditional moments.

The supporting space $Q$ can be found independently
of the conditional moments by means of analysis of 
a {\em moment matrix}.
Elements of a moment matrix are moments $M_\ell(B)$;
rows of moment matrix are indexed by response patterns $\ell$
having exactly one non-zero component;
columns of moment matrix are indexed by all possible response patterns.
The index of an element of the moment matrix lying in the intersection
of row $\ell'$ and column $\ell''$ is $\ell'+\ell''$.
However, addition of response patterns is possible only if
in every position either the first or second summand has a 0
(for example, $(1,0,0)+(0,2,1)=(1,2,2)$, but $(1,0,0)+(2,0,1)$
is undefined); thus, some elements of the moment matrix are
undefined.
The reason for some components being undefined is that we do not
have the possibility of performing a measurement on an individual 
multiple times independently, and since individuals are heterogeneous
(have different probabilities of outcomes of measurements),
we do not have multiple realizations of independent identically
distributed random variables.
In the example below, such components are shown by question marks.

Figure \ref{fig:MomMatr} gives an example of (a part of)
a moment matrix for the case $J=3$, $L_1=L_2=L_3=2$.
Columns in this matrix correspond to
$\ell=(000)$, $(100)$, $(200)$, $(010)$, $(020)$, $(001)$, $(002)$,
$(110)$;
other columns are not shown.

\begin{figure}[ht]   
\begin{equation*}
\begin{pmatrix}
M_{(100)}   &
    ?           &  ?            &
    M_{(110)}   &  M_{(120)}    &
    M_{(101)}   &  M_{(102)}    &
    ?           &  \cdots \phantom{\vdots}       \\
M_{(200)}   &
    ?           &  ?            &
    M_{(210)}   &  M_{(220)}    &
    M_{(201)}   &  M_{(202)}    &
    ?           &  \cdots \phantom{\vdots}       \\
M_{(010)}   &
    M_{(110)}   &  M_{(210)}    &
    ?           &  ?            &
    M_{(011)}   &  M_{(012)}    &
    ?           &  \cdots \phantom{\vdots}       \\
M_{(020)}   &
    M_{(120)}   &  M_{(220)}    &
    ?           &  ?            &
    M_{(021)}   &  M_{(022)}    &
    ?           &  \cdots \phantom{\vdots}       \\
M_{(001)}   &
    M_{(101)}   &  M_{(201)}    &
    M_{(011)}   &  M_{(021)}    &
    ?           &  ?            &
    M_{(111)}   &  \cdots \phantom{\vdots}       \\
M_{(002)}   &
    M_{(102)}   &  M_{(202)}    &
    M_{(012)}   &  M_{(022)}    &
    ?           &  ?            &
    M_{(112)}   &  \cdots \phantom{\vdots}      \\[-6pt]
&
\end{pmatrix}
\end{equation*}
\caption{\label{fig:MomMatr} Example of moment matrix}
\end{figure}   

For small $J$ (as in the example) one has large fraction
of undefined components in the moment matrix.
For large $J$ this fraction rapidly decreases.

The main fact with respect to the moment matrix is:

\begin{theorem}   
\label{th:MomMatr}
If a distribution is generated by $K$-dimensional LLS model
with supporting subspace $Q$, then
the moment matrix has a completion such that every one of its columns
belongs to $Q$.
\end{theorem}   

We show below that a completion of a moment matrix
is almost surely determined by its available part. Thus, the main
system of equations almost surely has a unique solution.

Here we need to clarify what we mean when we say ``uniqueness
of solution.''
The supporting subspace $Q$ and the identifiable properties
of mixing distribution on this subspace are defined uniquely.
One way to describe the subspace $Q$ is to present its basis,
and this can be done in infinitely many ways.
The {\em coordinate} expression of the mixing distribution
does depend on the choice of basis, and its characteristics
(like moments) also do depend on this choice.
Such dependencies are governed by tensor laws
\citep[][section 6.4]{Kovtun:2005}.

The existence of infinitely many bases does not mean that we have
infinitely many solutions;
rather, we can describe a {\em unique solution} by 
{\em infinitely many means}.
Various bases of $Q$ provide different points of view
on the same underlying picture.
The ability to choose a basis may benefit the applied researcher,
as it may help to present phenomenon under consideration more clearly.

The phase of finding the supporting subspace in LLS analysis
is tightly related to the principal component analysis
of the mixing distribution.
In fact, the submatrix of the moment matrix consisting of
columns from 2 to $|L|+1$ is (modulo incompleteness)
a shifted covariance matrix of the mixing distribution.
Theorem \ref{th:MomMatr} corresponds to the fact that
a multidimensional distribution is supported by $m$-dimensional
linear manifold if and only if the rank of covariance matrix is $m$.

Theorem \ref{th:MomMatr} also provides a method for determining
whether an LLS model exists for a particular dataset.
One has to find the largest computational rank \citep{Forsythe:1977}
of minors of the moment matrix containing no question marks;
if it is sufficiently smaller than $|L|$, LLS model exists
(the exact criterion of identifiability of LLS model is given
by theorem \ref{th:Ident}).

Now we are ready to describe a method for estimation of parameters
of LLS model.

First, the supporting subspace is estimated from the moment matrix.
The method of estimation is very similar to the one used in the
principal component analysis, adopted to handle incompleteness
of the moment matrix. The detailed description of the numerical
procedure is a subject of another article.

Second, a basis of the supporting subspace is chosen
and conditional moments $g^v_\ell$ are estimated by
(approximately) solving the main system of equations.
Note that moments can be found only for $\ell$ having sufficiently
many 0's (to guarantee that there are sufficiently many equations
(\ref{eq:MainEq})).
Moments for other $\ell$'s can be estimated as an average of
directly estimable moments (for example,
$g_{(\ell_1,\dots,\ell_J)} = \frac{1}{J} \left(
g_{(0,\ell_2,\dots,\ell_J)} + \dots + g_{(\ell_1,\dots,\ell_{J-1},0)}
\right)$).

This is a nonparametric approach, i.e. inference of properties of a
(mixing) distribution is made without any additional assumptions
about the structure of the distribution. 
If the nature of the applied problem justifies an assumption
that the mixing distribution belongs to some parametric family
(say, a mixing distribution is a Dirichlet distribution),
the parameters of such distribution may be easily estimated
by the moment method.

The mixed distribution in LLS analysis can be estimated in
style of empiric distribution, by letting the estimate of the
mixing distribution be concentrated in points $\bar{g}_\ell$
with weights $\bar{M}_\ell(B)$ (where bars mean estimates
of corresponding values).
It might be shown that this estimated distribution converges
to the true one when both size of the sample and number of measurements
tend to infinity, but the proof is outside the scope of the
present paper.


\section{Identifiability of LLS models} 
\label{sec:IdentifiabilityOfLLS}

Identifiability of a parameter of a model means that
the value of the parameter is uniquely determined by the
distribution of the observed variables \citep{Gabrielsen:1978}.
In our case, the observed distribution is given by moments $M_\ell(B)$.
Thus, identifiability means that the values of other parameters
(i.e., supporting subspace and conditional moments)
are uniquely determined by the values $M_\ell(B)$.

We start with the discussion of identifiability of supporting subspace
$Q$. The covariance matrix of the mixing distribution uniquely
defines $Q$ ($Q$ is spanned by the vector of expectations of
$\beta_{jl}$, which is the first column of the moment matrix,
and eigenvectors of covariance matrix corresponding to non-zero
eigenvalues). Thus, the supporting subspace is identifiable,
if the covariance matrix (which is incomplete for the same reasons
as the moment matrix)
can be uniquely restored from the available moments. 

\begin{lemma}   
\label{th:CovMatr}
Let $\mcC_m$ be a class of covariance matrices of $n$-dimensional
distributions of rank $m$.
Let, for arbitrary $A \in \mcC_m$, $\tilde{A}$ denote the matrix
$A$ with missing diagonal blocks of maximal size $p$.
Let also assume that inequality $2m+2p-1 \le n$ holds.

Then, for arbitrary $A,B \in \mcC_m$,
the equality $\tilde{A} = \tilde{B}$ almost surely implies $A=B$.
\end{lemma}   

\begin{proof}[Outline of the proof]
Figure \ref{fig:CovMatr} demonstrates how an incomplete
covariance matrix can be restored. 
If the minor $(c^i_j)_{ij}$ is nondegenerate, there exist a unique
linear combination of columns $c_1,\dots,c_m$ that yields column $b$;
let $b_j = \sum_i \gamma_i c^i_j$ for all $j$.
Then, as the rank of the whole matrix is $m$, the element in the top
left corner of the matrix must be $\sum_i \gamma_i a_i$.
Other elements denoted by question marks may be restored by
applying similar procedure.

\begin{figure}[ht]   
\begin{equation*}
\begin{matrix}
\hfill \le p \text{ rows } \left\{ 
    \vphantom{\begin{matrix}\vdots\\\vdots\end{matrix}} 
        \right. \\
\vphantom{\begin{matrix} b_1 \\ \vdots \\ b_m \end{matrix}} \\
\text{some rows } \left\{ \vphantom{\vdots} \right. \\
\vphantom{ \boxed{\begin{matrix} ? \\ \vdots \end{matrix}} }  \\
\vphantom{\begin{matrix} \ddots \\ \boxed{?} \end{matrix}}
\end{matrix}
\begin{pmatrix}
\boxed{ \begin{matrix} ? & \dots \\ \vdots & \ddots \end{matrix} } &
& 
\begin{matrix} & a_1 & \dots & a_m \\ & & \vdots & \end{matrix} \\
\begin{matrix} b_1 & \\ \vdots & \\ b_m & \end{matrix} &
\begin{matrix} & \\ & \ddots \\ & \end{matrix} &
\begin{matrix} & c^1_1 & \dots & c^1_m \\ & \hdotsfor[2]{3} \\
               & c^m_1 & \dots & c^m_m \end{matrix}  \\[+25pt]
\hdotsfor[3]{3} \\[+5pt]
& & \begin{matrix} 
\boxed{ \begin{matrix} ? & \dots \\ \vdots & \ddots \end{matrix} } &
\phantom{\boxed{ \begin{matrix} ? & \dots \\ \vdots & \ddots \end{matrix} } }
\end{matrix} \\
& & \begin{matrix}
\phantom{\boxed{ \begin{matrix} ? & \dots \\ \vdots & \ddots \end{matrix} } }
& \begin{matrix} \ddots & \\ & \boxed{?} \end{matrix}
\end{matrix}
\end{pmatrix}
\begin{matrix}
\vphantom{ \boxed{\begin{matrix} ? \\ \vdots \end{matrix}} } \\
\left. 
    \vphantom{\begin{matrix}c^1_1\\ \hdotsfor{1} \\ c^1_1 \end{matrix}}
        \right\}~m \text{ rows} \hfill \\
\vphantom{\}} \\[+15pt]
\left. 
    \vphantom{ \boxed{\begin{matrix} ? \\ \vdots \end{matrix}} }
        \right\}~\le p \text{ rows} \hfill \\[+15pt]
\left. 
    \vphantom{\begin{matrix} \ddots \\ \boxed{?} \end{matrix}}
        \right\}~\le (m-1) \text{ rows}
\end{matrix}
\end{equation*}
\caption{\label{fig:CovMatr}Restoration of a covariance matrix}
\end{figure}   

The picture also illustrates the necessity of condition $2m+2p-1 \le n$.

Thus, to complete a proof of the theorem, it is sufficient to show
that all minors of a covariance matrix are almost surely nondegenerate.

Every covariance matrix $A$ of rank $m$ is nonnegative definite and
symmetric; thus, it can be represented in form $A = O^T D O$,
where $O$ is an orthogonal matrix and $D$ is a diagonal matrix with
exactly $m$ nonzero elements. Further, almost every orthogonal matrix 
may be represented as $O = (I+V)(I-V)^{-1}$, where $I$ is the unit
matrix and $V$ is a skew-symmetric matrix 
\citep[Cayley parametrization; see][IV.6]{Satake:1975}.
Thus, elements of covariance matrix can be represented as
ratios of polynomials of $\frac{n(n-1)}{2}+m$ variables
($\frac{n(n-1)}{2}$ elements of skew-symmetric matrix and $m$
elements of matrix $D$).
Consequently, all minors of order $m$ are ratios of polynomials
of these variables, and they are not identically 0 (as it is possible
to give an example of covariance matrix of rank $m$ with all minors
of order $m$ being non-degenerate).
But a set of 0's of a polynomial has measure 0, q.e.d.
\end{proof}

\begin{corollary}
\label{cr:CovMatr}
Let $\mfM$ be a family of mixing distributions supported by
$K$-dimen\-sio\-nal linear subspaces and let $\nu$ be a Borel measure
on $\mfM$. Let $\phi : \mu \mapsto C_\mu$ be a mapping of mixing
distributions to their covariance matrices.
Suppose the image measure of $\nu$ under mapping $\phi$ is
absolutely continuous with respect to Lebesgue measure.
Then for $\nu$-almost all mixing distributions their supporting
subspace is identifiable.
\end{corollary}

\begin{remark}
For a finite-dimensional Euclidean space, there is a ``standard''
measure (Lebesgue measure), with respect to which ``almost surely''
statements are usually made. Unfortunately, there is no such
``standard'' measure in the infinite-dimensional case, and particulary
there is no ``standard'' measure on the space of mixing distributions.
However, one can introduce a notion of ``nowhere degenerated'' measure
and show that every ``nowhere degenerated'' measure satisfies
conditions of corollary \ref{cr:CovMatr}.
Thus, one can say that the supporting subspace of the mixing
distribution is almost surely identifiable with respect to any
``nowhere degenerated'' measure on the space of mixing distributions.
\end{remark}

As lemma \ref{th:CovMatr} shows, the supporting subspace of
the mixing distribution is (almost surely) uniquely defined 
by moments $M_\ell(B)$.
This linear subspace can be described by infinitely many bases.
When a basis of the supporting subspace is chosen, the main
system of equations becomes a linear system with respect
to conditional moments $g^v_\ell$. 
Moreover, it breaks apart into small subsystems that can be solved
separately. For example, conditional expectations $g_{\ell k}$
satisfy equations (\ref{eq:MainEqExp}).
If $\ell$ contains zeros at places $j_1,\dots,j_p$, one
has $l(\ell)=L_{j_1}+\dots+L_{j_p}$ equations for 
$g_{\ell 1},\dots,g_{\ell k}$, from which $l(\ell)-p$ are
independent.
Thus, $g_{\ell k}$ may be uniquely determined from the system,
if $k \le l(\ell)-p$.
Other conditional moments can be uniquely calculated from the 
system under similar conditions.

Summarizing, we obtain

\begin{theorem}
\label{th:Ident}
If the observed distribution is generated by a $K$-dimensional
LLS model with $K \le \frac{|L|-J}{2}-\max L_j + \frac{5}{2}$.
Then the LLS model is almost surely identifiable.
\end{theorem}


\section{Consistency of LLS estimators} 
\label{sec:ConsistencyOfLLS}

The consistency of LLS estimators is almost a straightforward corollary
of the well-known statistical fact that frequencies are consistent
and efficient estimators of probabilities.

The supporting subspace is estimated as a $K$-dimensional subspace
closest to columns of the frequency matrix (more precisely, closest
to the subspaces spanned by incomplete columns).
This estimate continuously depends on the elements of the frequency
matrix; thus, it converges to the true supporting subspace
when elements of the frequency matrix converge to the true moments.

Similarly, estimators for conditional moments $g^v_\ell$ are
(approximate) solutions of a linear system with coefficients 
depending on frequencies. Again, these estimators continuously depend
on frequencies and converge to the true conditional moments when
frequencies converge to the true moments.

The consistency of LLS estimators may be formulated as follows.
Suppose we have a distribution generated by LLS model with
mixing distribution supported by subspace $Q$ and having
conditional moments $g^v_\ell$.
Then estimators $\bar{Q}$ and $\bar{g}^v_\ell$, obtained by the
procedure described above, converge to $Q$ and $g^v_\ell$,
respectively, when the size of a sample tends to infinity.
Thus, we have:

\begin{theorem}
If LLS model is identifiable, it is consistently estimable.
\end{theorem}


\section{Simulation studies} 
\label{sec:SimulationStudies}

We have developed a prototype of the algorithm for estimation of
LLS parameters and performed preliminary experiments with it.
The implementation of the algorithm follows the ideas described above,
though differring in detail needed to provide computational stability.

\begin{figure}[ht]
\centering
\scalebox{0.3}{\includegraphics{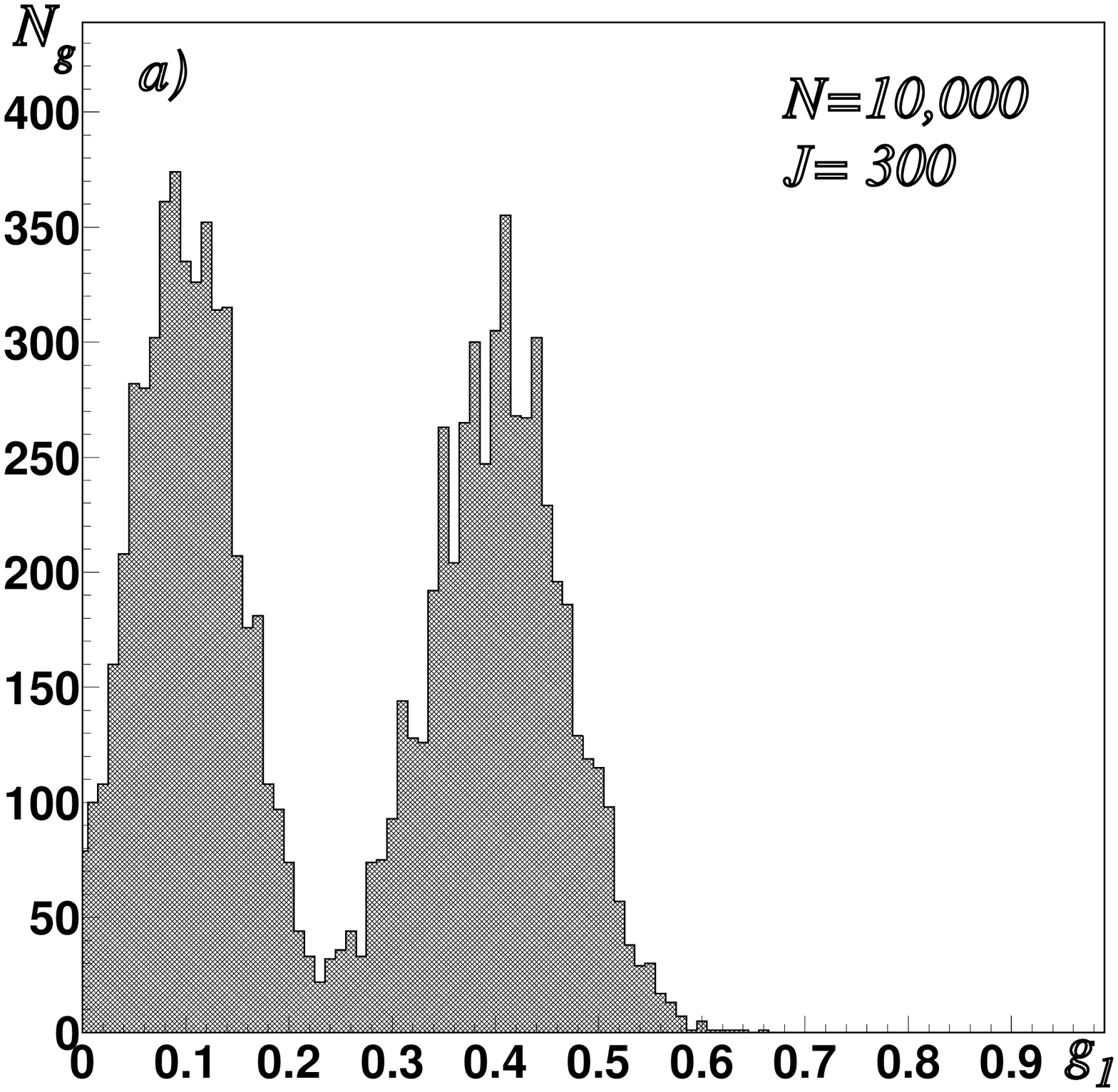}}
\hfill
\scalebox{0.3}{\includegraphics{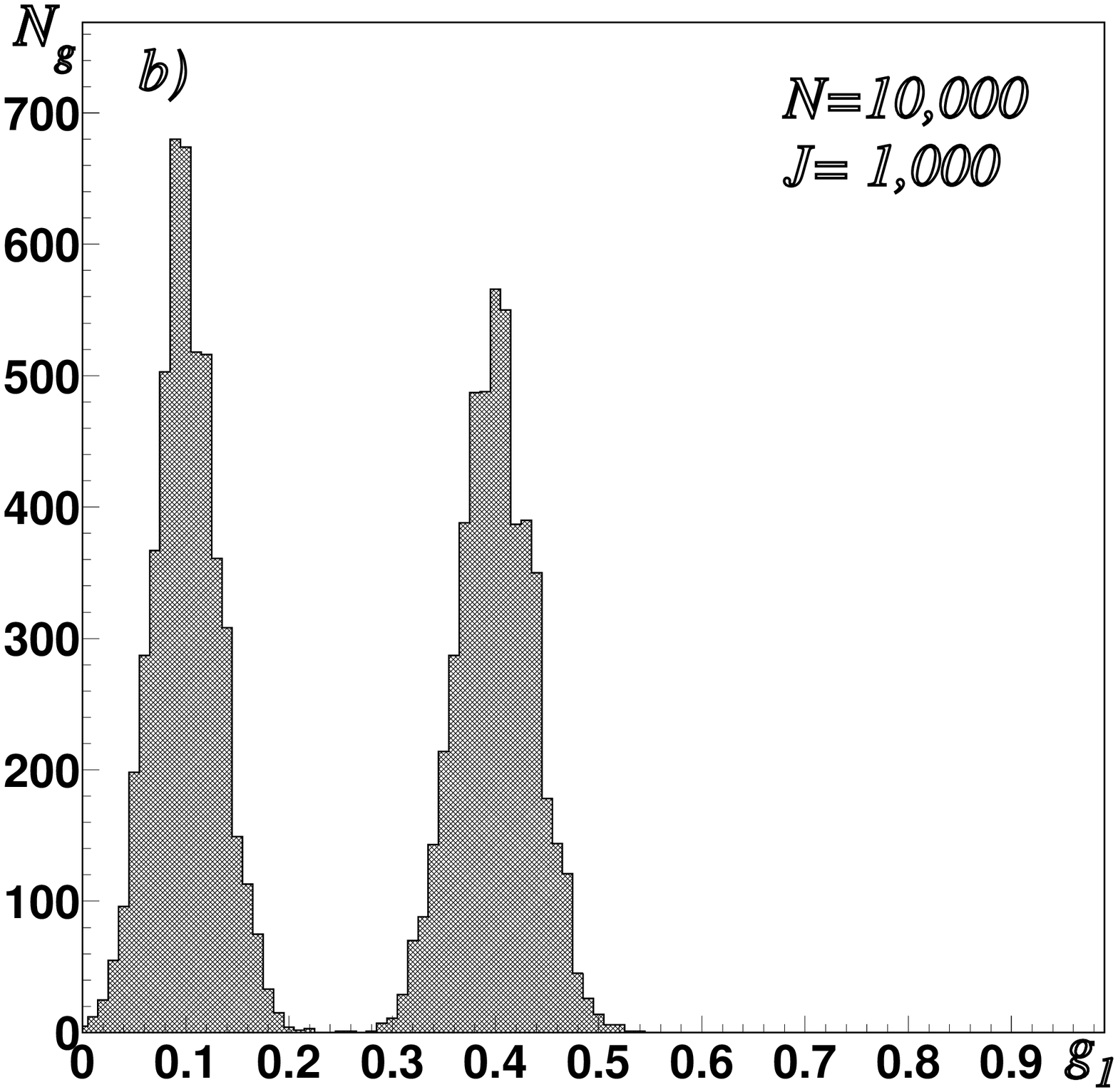}}
\\ \medskip
\scalebox{0.3}{\includegraphics{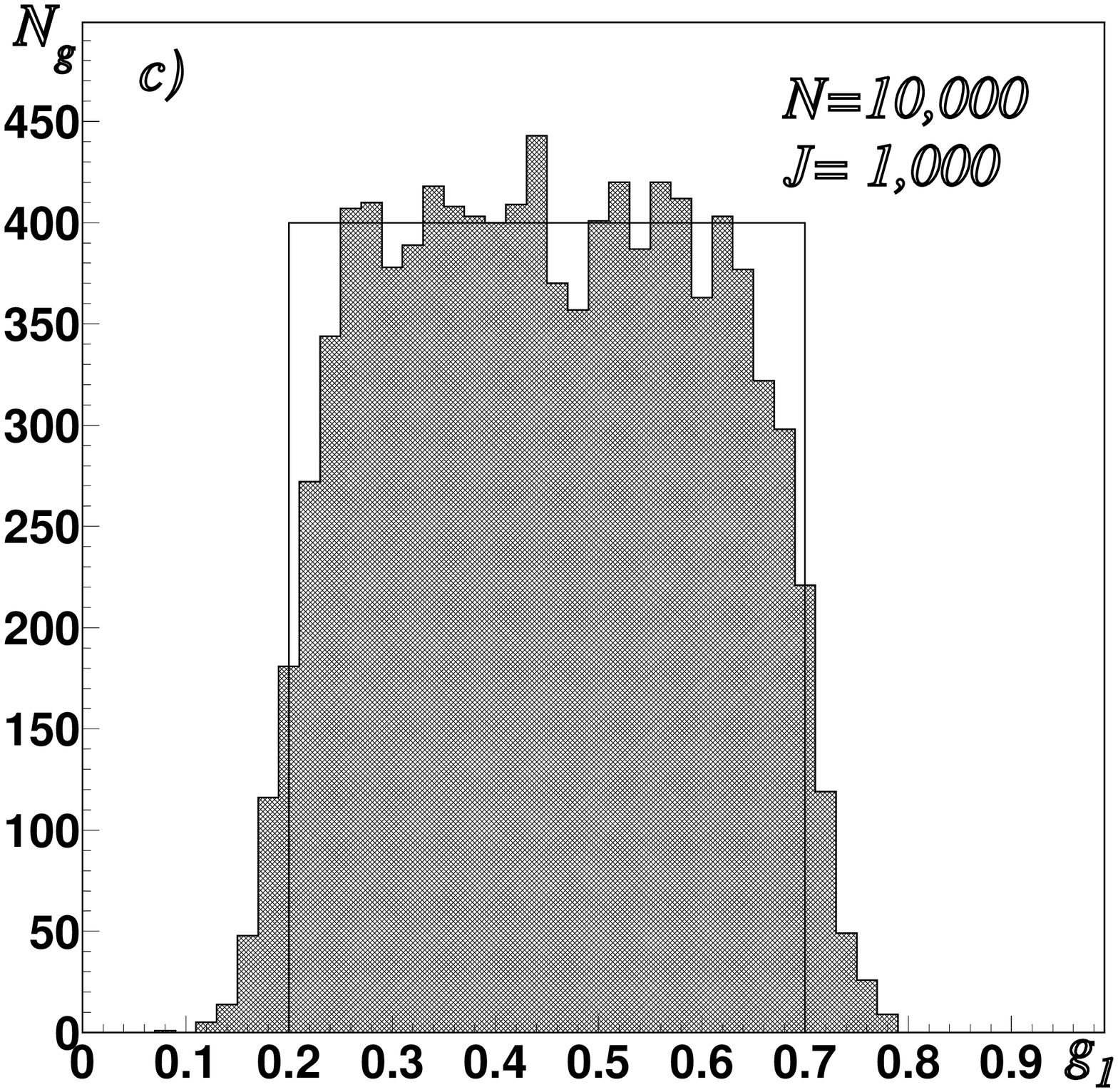}}
\hfill
\scalebox{0.3}{\includegraphics{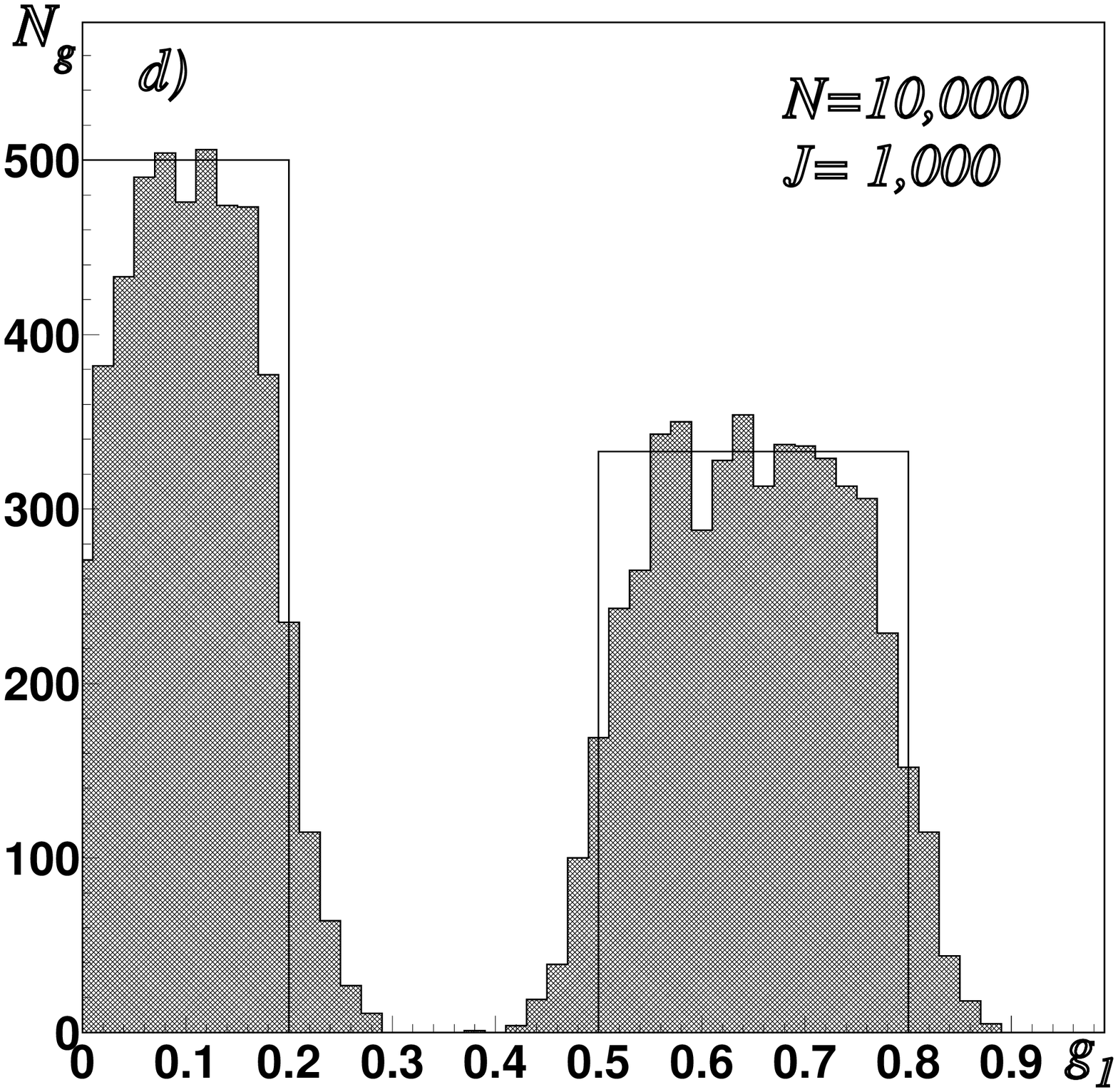}}
\caption{
    \label{fig:RestoredDistr}
     Restored mixing distribution
     (see text for further explanations).
}
\end{figure}

The first experiments with the algorithm gave encouraging results.
For illustrative purposes, we choose 2-dimensional LLS model.
As LLS uses homogeneous coordinates $g=(g_1,g_2)$ and $g_2 = 1-g_1$,
this means that the mixing distribution can be thought of as a
distribution over interval $g_1 \in [0,1]$.
The results of four experiments are presented in figures
\ref{fig:RestoredDistr}a--d.
All experiments were organized as follows.

We randomly generated 2 basis vectors of the supporting subspace.
Figure \ref{fig:RestoredDistr}a is a case with 300 binary questions
(i.e., $J=300$, $L_j=2$ for all $j$, $|L|=600$);
figures \ref{fig:RestoredDistr}b--d are cases with 1000 binary questions
(i.e., $J=1000$, $L_j=2$ for all $j$, $|L|=2000$).
Then, we choose a mixing distribution.
In figures \ref{fig:RestoredDistr}a,b the mixing distribution
is concentrated at two points, $0\text{.}1$ and $0\text{.}4$,
in figure \ref{fig:RestoredDistr}c the mixing distribution
is uniform over subinterval $[0\text{.}2,0\text{.}7]$,
and in figure \ref{fig:RestoredDistr}d it is uniform over two
subintervals, $[0,0\text{.}2]$ and $[0\text{.}5,0\text{.}8]$.
Then we generated a sample of 10,000 individuals by
randomly choosing a point in the supporting subspace
in accordance with the mixing distribution
and generating responses with probabilities defined by the selected
point (using equation (\ref{eq:BetaG})).

The set of responses was used as input to the algorithm.
The algorithm estimated the supporting subspace, conditional
moments of mixing distributions, and then estimated the mixing
distribution itself. Histograms of restored distribution
are shown in figures \ref{fig:RestoredDistr}a--d; a solid line
in figures \ref{fig:RestoredDistr}c--d shows the true mixing
distribution used to generate samples.

Experiments with different choices of supporting subspace
and other randomly generated samples give similar results.

Figures \ref{fig:RestoredDistr}c--d demonstrate a good quality
of restoration of the mixing distribution under various conditions.
Figures \ref{fig:RestoredDistr}a--b show that the
precision of restoration of the mixing distribution increases with the
increase of the number of variables.

It is interesting to compare the results of LLS model with the
results of latent class model (LCM).
In the cases \ref{fig:RestoredDistr}a and \ref{fig:RestoredDistr}b,
LCM would restore the same picture as LLS does, and thus, it can be
an alternative to LLS model. However, in the cases
\ref{fig:RestoredDistr}c and \ref{fig:RestoredDistr}d,
LCM may be used only as approximation, and it might be shown
that LCM may involve approximately 1,000 (number of measurements)
latent classes in these cases.


\section{Discussion} 
\label{sec:Discussion}

In the present paper we described a new class of models for analyzing
high-dimensional categorical data, which belongs to a family
of latent structure models.
We established conditions for identifiability of models
and consistency of parameter estimators.
The essence of our approach is in the consideration of a space
of independent distributions (which is a space of distributions
being mixed in the case of latent structure analysis)
as a linear space.
This allows us, first, to formulate model assumptions in the language
of linear algebra, and second, to reduce a method for estimating
model parameters to a sequence of linear algebra problems.
The very modest identifiability conditions (theorem \ref{th:Ident})
allow application of these models to a wide range of practical datasets.

This linear-algebra approach allows us to clarify relationship
between various branches of latent structure analysis.
Consider, for example, relation between LLS models and
latent class models (LCM). 
In geometric language, latent classes are points in the space
of independent distributions.
If an LCM with classes $c_1,\dots,c_m$ exists for a particular dataset,
then an LLS model also exists, and its supporting subspace is
the linear subspace spanned by vectors $c_1,\dots,c_m$.
Thus, dimensionality of LLS model never exceeds the number of classes
in LCM. These numbers are equal if and only if LCM classes
are points in general position
(i.e., vectors $c_1,\dots,c_m$ do not belong to a subspace
of dimensionality smaller than $m$).
If LCM classes are not in general position, however, the dimensionality
of LLS model may be significantly smaller. For example,
it is possible to construct a mixing distribution such that
(a) it is supported by a line (i.e., dimensionality of LLS model
is 2); (b) there exists LCM with $J$ (number of variables) classes;
(c) there is no LCM with smaller number of classes.
(A rigorous proof of the last fact will be given in another paper.)
On the other hand, LLS can be used to evaluate applicability of
LCM: if the mixing distribution in LLS model has pronounced
modality, then an LCM is more likely to exist (with the number 
of classes equal to number of modes).

Maybe, the most important question regarding any kind of model
is its interpretation. The interpretation heavily depends on
application domain, so we are able to give here only very general
guidelines.
If the application domain supports an assumption that
individuals in a population may be described by points in a
state space and probabilities of outcomes of measurements
depend on individual coordinates in the state space,
this state space can be recovered by LLS analysis,
and coordinates of an individual in the state space can be estimated
from the outcomes of measurements.
However, the ``physical meaning'' of the state space,
what does it mean ``to be in a particular region of the state space'',
etc. may be discussed only in terms of the application domain.

There is one interesting property of LLS models,
which can be characterized as a {\em partial identifiability}.
Our method gives consistent estimates for supporting subspace
of the mixing distribution and for conditional moments
$g^v_\ell$ of maximal order $v$ satisfying $|v|=v_1+\dots+v_K \le J$.
This is not, however, a limitation of the model;
rather, it is limitation of the problem itself:
if two mixing distributions are supported by the same subspace
and have the same conditional moments of order $|v| \le J$,
they will produce the same observed moments $M_\ell$;
thus, these two mixing distributions are indistinguishable based
on available data.
On the other hand, two mixing distributions, which have the
same moments of order up to $J$ cannot be {\em significantly
different}: it might be shown that distance between them
converges to 0 when $J$ tends to infinity.
This means that the recovered knowledge about mixing distribution
can be made more and more precise by increasing
number of measurements.
This fact is well recognized in practice; for example,
a mathematical test based on multiple-choice questions
would include several questions regarding, say, the addition of fractions
to judge student performance on this topic.

The above problem may also be considered from another side.
As it was mentioned in the end of section
\ref{sec:EstimationOfLLSModel}, the mixing distribution can be
estimated in style of empirical distribution,
by letting the estimate of the mixing distribution be concentrated
in points $\bar{g}_\ell$ with weights $\bar{M}_\ell(B)$.
One can ask how this estimate relates to the true mixing
distribution?
The answer is: the estimate of the mixing distribution
converges to the true one when both size of a sample and number
of measurements tend to infinity.
This fact may be considered as an analogue of the
Glivenko-Cantelli theorem.
The fact that estimate of the individual position in the state space
becomes more and more precise with the increase of the number
of measurements is an analogue of the Bernoulli's law
of large numbers.
The fact that one needs more and more measurements
performed on each individual to increase precision of restoration
of the mixing distribution does not diminish the usefulness
of LLS analysis.
It is well recognized that to achieve a required precision in
statistical inference one needs to perform sufficiently many
measurements. The difference here is in that one needs not
only to repeat the same measurement on different individuals,
but also to perform sufficiently many measurements on each individual.
The proof of the above convergence and estimation of the rate
of the convergence is subject of forthcoming papers.


\bibliographystyle{apalike}
\bibliography{math,probability,LSA,GoM,mkovtun}

\end{document}